\input amstex
\documentstyle{jmhppt}
\input epsf
\magnification=1200
\hsize=13cm
\vsize=19cm
\def\h{Hammersley}
\def\msection#1{\bigskip\centerline{\uppercase{#1}}\smallskip\noindent\is}
\def\ssection#1{\bigskip\centerline{\it#1}\smallskip\noindent\is}
\def\normalparindent{12pt}
\parindent=\normalparindent

\catcode`\@=11
 \def\logo@{}
\catcode`\@=13

\font\headfont=cmr8

\footline={\hfil}
\headline={\ifnum\pageno<3{\hfil}\else{\headfont John Michael Hammersley (1920--2004) \hfil \folio}\fi}

\def\is{\ignorespaces}
\long\def\quote#1\endquote{{{\it#1}}}
\def\ii{\itemitem}
\def\iii#1#2{\parindent=2\parindent\ii{{#1}}{#2}\par
             \parindent=\normalparindent}
\input jmhlab
\def\cite#1{[\citeto{#1}]}
\def\citem#1{\citeto{#1}}
\def\pc{p_{\text{\rm c}}}
\def\O{\text{\rm O}}

\def\jmh{Hammersley, J.\ M.\ignorespaces}
\def\nref#1\endref{}

%\topinsert
\line{\hfil}
\vskip2cm
\centerline{\epsfxsize=11cm \epsfbox{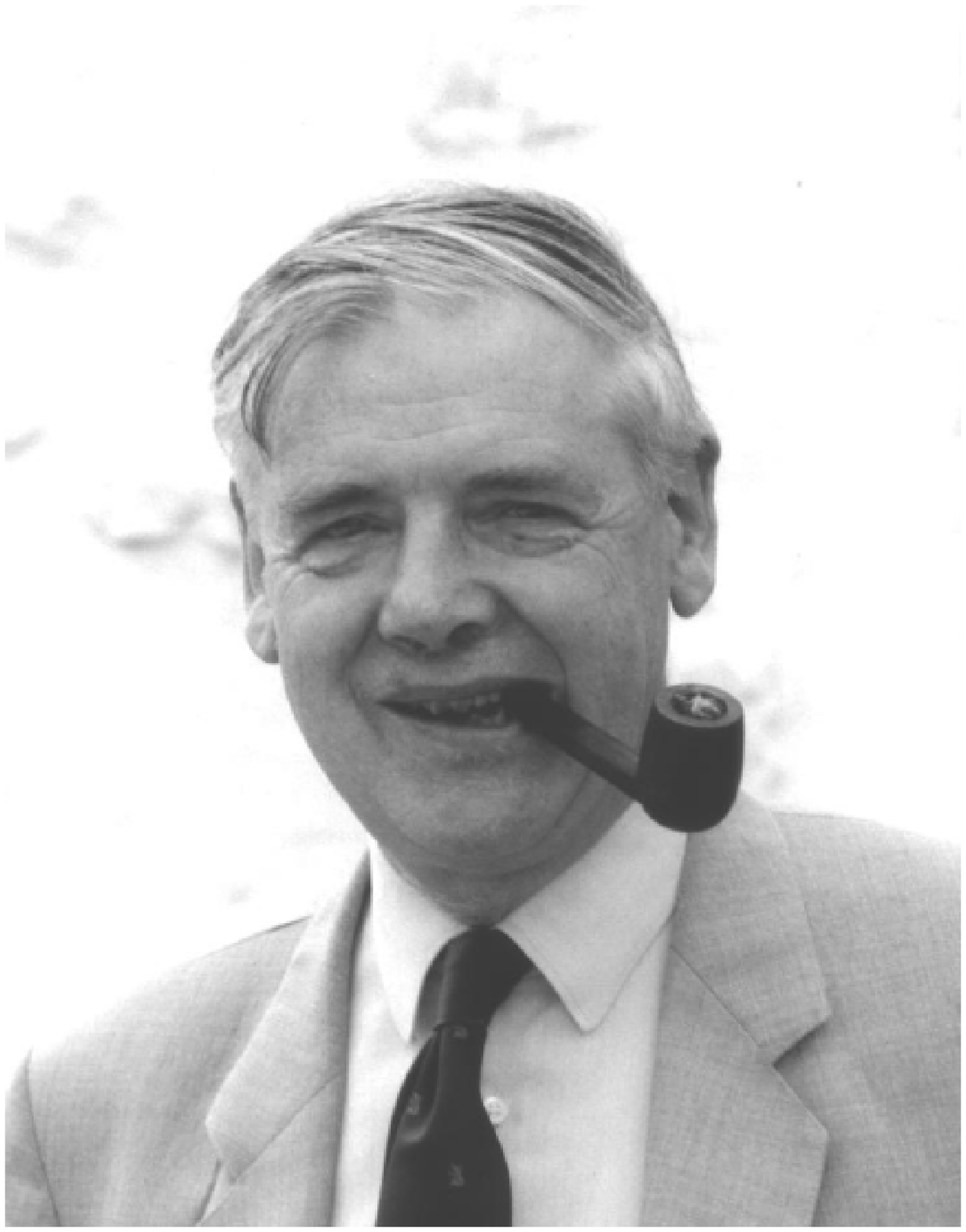}} 
%\centerline{\vbox{(photo here)}}
\botcaption{John Michael Hammersley}
\endcaption
%\endinsert
\vfill\eject

\centerline{JOHN MICHAEL HAMMERSLEY}
\centerline{21 March 1920 --- 2 May 2004}
\centerline{Elected FRS 1976}
\medskip
\centerline{By Geoffrey Grimmett and Dominic Welsh}
\medskip
\centerline{\it Centre for Mathematical Sciences, 
University of Cambridge, Cambridge CB3 0WB}
\centerline{\it Merton College, Oxford OX1 4JD}

\bigskip

\noindent
John \h\  was a pioneer amongst mathematicians who defied classification
as pure or applied; when introduced to guests at Trinity College, Oxford,
he would say he did \lq\lq difficult sums". 
He believed passionately in the importance of mathematics
with strong links to real-life situations,
and in a system of mathematical education in which the solution of problems
takes precedence over the generation of theory.
He will be remembered for his
work on percolation theory, subadditive stochastic
processes, self-avoiding walks, and Monte Carlo
methods, and, by those who knew him, for his intellectual
integrity and his ability to inspire and to challenge.
Quite apart from his extensive research achievements,
for which he earned a reputation as
an outstanding problem-solver, he was a leader in the movement of the 1950s
and 1960s to re-think the content of school
mathematics syllabuses.

\msection{Family background}
John \h\ was born to a couple with strong international connections.
His mother,
Marguerite (n\'ee Whitehead), was 
born on 29 June 1889 in Moscow, where her father Thomas was
engaged in the export and sale
of cotton-spinning and other textile machinery from Lancashire.
At the age of 14, she was sent to boarding school in England,
thus escaping the difficulties and deprivations
faced by
her brothers, and documented in \cite{new56}, 
as a consequence of the Revolution of 1917, when
the Bolsheviks declared all foreign assets to be owned by the Russian people.
Their property was seized, and their families retreated to London
via Murmansk in 1918. Early on 1 January 1920, John's uncle George was hauled
out of bed by the secret police (the Cheka) and interrogated
over a period of three weeks in the Lubianka, sleeping on a bare concrete floor
at sub-zero temperatures. George's brother, Alfred,
managed to extricate George from the labour camp to which
he had been moved, on the grounds that
he was about to die. He survived, however, and he and Alfred
caught a train that same day to the Finnish border.

John's father, Guy Hugh, was born on 5 March
1883, the second son of a fashionable London
gynaecologist who, when Guy was 14, 
collapsed and died in the prime of life, leaving
his family in straitened circumstances. 
Guy had to leave school, and he took a 
job as an office boy at the London office of the United States
Steel Corporation.
By the time of John's birth in 1920, 
Guy had worked his way up to be in charge of the branch office 
in Glasgow. There were ups and downs in his career, occasioned by
times of retrenchment and recession in the United States.
Guy and Marguerite moved back to England, and
he was made redundant around 1925.  He found work as the London manager
for the Youngstown Steel and Tube Company, and later as
European manager for the Bethlehem Steel Company following the 
Depression in the USA.

Marguerite and Guy were married in 1914, and their only son
to survive childbirth, John Michael, was born on 21 March 1920.

\msection{Education}
The following extracts from some autobiographical
notes present an interesting account of John's life pre-Oxford,
as well as insight into his character.

\quote
I attended a kindergarten called the Waterside School in Bishops 
Stortford from 1925 to 1929. It was run by a headmistress, Miss 
Blandford, and it gave me an excellent start in reading and writing 
and arithmetic. In my last year, her father, Mr Blandford, 
gave me an introduction to Latin and algebra.

In 1929 I was sent as a boarder to Bembridge School on the Isle
of Wight. This was a school with progressive ideas about
teaching arts and
crafts and carpentry but little emphasis on anything academic: 
after a couple of terms at Bembridge,
my parents were dissatisfied with what I was being taught and 
I was sent instead to a more conventional
preparatory school, Stratton Park near Bletchley,
where I remained from 1930 to 1934.

The man who taught mathematics at Stratton Park, Mr Pilliner,
almost put me off the subject by asking me how many blue beans made five.
When I failed to answer the conundrum, he said the
answer was $5$ and I was a fool: but I had already dismissed this 
as too obvious to be correct (and in retrospect, the correct answer is
probably something like $5[\text{\rm blue beans}]^{-1}$).
However,
my mathematical fortunes were saved shortly after this incident by the arrival
at Stratton Park of another teacher
of mathematics, Gerald Meister. He had been a housemaster
at Sedbergh School, where there was a convention that
housemasterships could only persist for 15 years. When his 15-year
stint was complete, he decided to try his hand at
preparatory school teaching and took up residence at Stratton Park 
and remained there for a couple of years, after which 
he taught at Wellington College
and next at the Dragon School in Oxford.

During his time at Stratton Park he gave me a solid education in mathematics
and a liking for the subject. This covered plenty of Euclidean
geometry (including such topics as the nine-point circle) and algebra
(Newton's identities for roots of polynomials) and trigonometry (identities
governing angles of a triangle, circumcircle, incircle, etc),
but no calculus. Due to his help, I got a scholarship
to Sedbergh.

I was at Sedbergh School from 1934 to 1939.
There it was traditional in those days for the brighter boys 
to be shoved on the classical side, and in my first year
I was put in the Classical Fifth form, where I completed 
the School Certificate
in classics (the equivalent of four O-levels today) 
and then at the end of my first year into
the Lower Sixth Classical. However, Latin and Greek did not
interest me, and after one term in the Lower
Sixth Classical I was allowed to migrate to the Upper Sixth
Modern to learn some science. 
I had some excellent teaching
in physics from Len Taylor, and in chemistry from Charles \hbox{\rm[sic]}
Mawby\footnote{N.\ James Mawby}.
My mathematics master was Sydney Adams (subsequently headmaster
of Bancrofts School). His knowledge of mathematics
was very sound, but did not extend much beyond what was appropriate
to schoolteaching: I recall being puzzled that a continuous function might be
non-differentiable everywhere; and although he was able to
confirm this, he could not exhibit a specific example for me.
I passed Higher Certificate (the equivalent of A-level today)
in mathematics, physics, and chemistry in the summer of 1937, but
I did not gain a distinction in mathematics. I sat the
scholarship examination for Emmanuel College, Cambridge, in
December 1937, and
also for New College, Oxford, in March 1938, without
success in both cases. However, I was awarded a Minor Scholarship to
Emmanuel College at a second attempt in December 1938.
\endquote

\msection{Cambridge}
\quote
I went to Cambridge as an undergraduate in 1939. 
The war had just started,
and many undergraduates including myself presented ourselves to enlist at the
Senate House which served as a recruiting station in Cambridge.
At least as far as this recruiting station was concerned, there was not
much evidence at that time of making
wartime use of people with potential scientific qualification.
After a brief medical check-up, I found myself
in front of a trestle table opposite a don, disguised
in the uniform of a sergeant, and the following conversation ensued.
\iii{Sergeant:}{Do you want to join the navy, the army, or the air force?}
\iii{Me:}{I suppose it should be the army --- I was in the OTC at school.}
\iii{Sergeant:}{Which regiment do you have in mind?}
\iii{Me:}{I have no idea. I have just started to read
mathematics here in Cambridge: is there any use for mathematics
in the army?}
\iii{Sergeant:}{No, there is no
use for mathematics in this war and in any case you are only an undergraduate.
The services have taken just three professional mathematicians
from Cambridge, one for the navy to tell them about underwater
explosions, one for the air force to explain stellar navigation,
and I was the third. My mathematical job is
to add up the daily totals of recruits for the navy, the army
and the air force respectively.}

\noindent
I wonder who the `sergeant' was, maybe a number theorist.
Of course, he was wrong\footnote{He was more or less in agreement
with G.\ H.\ Hardy, \cite{new57}, who felt it plain
that ``the real mathematics has no {\it direct\/} utility in war",
but, when asking ``does mathematics ``do good" in war?",
found it probable that technical skill keeps young
mathematicians from the front, thereby saving their lives.}
 about the wartime uses of science, 
including mathematics,
and about the number of scientists and mathematicians 
recruited from Cambridge, but I did
not know about that until much later. In the meantime, waiting
until I was eventually called up, I hung around in Cambridge pretty
idly. I remember tutorials from Stonely,
who taught me how to express $\nabla^2\phi$ in spherical polar coordinates
but not much else; and also tutorials from P.\ W.\ Woods,
whose favourite subject was the twisted cubic. Pupils would
strive to keep him off the twisted cubic for as long as possible by
asking him questions on other
bits of pure mathematics, but once he was locked on the twisted cubic
after the first ten minutes of a tutorial, the rest of the tutorial was a
foregone conclusion. I was lucky to get
a Third Class\footnote{Of 33 candidates for the 
Mathematics Preliminary Examinations in 1940,
11 were placed in the First Class, 15 in the Second, and 7 in the Third.}
in the preliminary examinations in mathematics
in the Easter Term in 1940, before being called up for military service in the
Royal Artillery.
\endquote

\msection{Wartime Service}
\quote
Despite the assertion by the recruiting `sergeant' in the Senate House 
in Cambridge that mathematics was of no military interest in 
wartime, I did later find uses for it when serving in the Royal 
Artillery in connection with anti-aircraft gunnery. 
An aircraft is a high-speed moving target, whose flight path is 
detected and followed by radar. To hit a target one needs to predict 
how far the aircraft will have moved in the lapse of time between the 
gun being fired and the shell reaching it. This calculation was 
performed by a piece of computing hardware called a predictor. 
There were two  sorts of staff officers who were expected to have 
an enhanced technical knowledge of anti-aircraft equipment: 
they were respectively called Instructors in Gunnery (I.G.s) and 
Instructors in Fire Control (I.F.C.s) Both the I.G.\ and the I.F.C.\  
had technical expertise in the three components (radar, predictor, gun) 
of this linkage; but their particular provinces overlapped in the sense 
that the I.G.s specialized in the gun--predictor pair, while the I.F.C.s 
specialized in the radar--predictor pair. The School of Anti-Aircraft 
Artillery (S.A.A.A.) was situated on the Pembrokeshire coast at 
Manorbier; and the Trials Wing of the S.A.A.A.\ was at Lydstep 
about a mile to the east of Manorbier. The function of the Trials 
Wing was to carry out research on the performance of various pieces 
of anti-aircraft equipment, both  existing equipment and equipment 
proposed for future use, and to report thereon to the war office and 
Ministry of Supply. At the Trials Wing there were three I.G.s and 
two I.F.C.s; and in 1942 I became one of the I.F.C.s, 
remaining there until the end of the war.

Before that however I was called up for military service in the late 
summer of 1940, first as a gunner and next as a lance-bombardier 
at a training camp at Arborfield until being sent to an officer 
training cadet unit at Shrivenham. I was commissioned as a second 
lieutenant in the spring of 1941 and posted to an anti-aircraft 
gun site defending an armament factory near Warsham.  
At Shrivenham I had been told about the existence of radar; 
and the Warsham gun site had an early piece of radar equipment 
which operated with a wavelength of a few metres. Its performance 
in measuring the distance to a target was reasonable; but its 
accuracy in measuring the direction to the target was pretty 
indifferent, relying on interference effects between various dipole 
aerials receiving signals both directly and also reflected from a 
large horizontal mat of wire mesh. At any rate it represented the 
current state of the art at that time; and it interested me 
considerably. Wanting to learn more about the potentialities of 
radar, I took the rather unusual step of telephoning divisional 
headquarters and as a result was selected to train to become an I.F.C.

This training began with a six weeks course on basic wireless 
technology at the Regent Street Polytechnic, followed by a longer 
and more specialized course on radar at Watchet in Somerset. 
At Watchet they had a radar with a ten centimeter wavelength, 
which at that time had not come into general service for 
anti-aircraft gunnery. There I learnt about the properties of 
magnetrons and wave guides. On passing out of Watchet as a qualified 
I.F.C., which carried the automatic rank of captain, 
I was posted first to an establishment at Oswestry which trained 
operators of radar equipment, and next to anti-aircraft brigade 
headquarters in the Orkneys where I was responsible for the radar 
installations of the gun sites defending Scapa Flow. 
Finally in 1942 I was transferred to the Trials Wing at Lydstep.

Amongst the personnel at Trials Wing there was a team of about 40 girls 
who carried out the computations necessary for analyzing the 
performance of the anti-aircraft equipment, and I was responsible 
for directing their calculations. One of their jobs consisted in 
operating the kinetheodolites for tracking a target. The 
kinetheodolites were a pair of synchronized telescopic cameras at 
each end of a base line about a couple of miles long, which could 
give simultaneous readings of the respective angles to a target 
(either an aircraft or a radar sleeve towed behind an aircraft). 
From the resulting data it was possible to compute fairly accurate 
positions of the target and how these positions  
depended upon time as the target moved along its flight path. 
In practice it was just an ugly piece of three-dimensional trigonometry; 
and when I first arrived at Lydstep it was done with pencil and paper 
with the aid of a 7-figure tables of trigonometric functions, 
in accordance with traditions of military surveyors. 
But while surveyors may conceivably be interested in determining 
a position to the nearest fraction of an inch, it was nonsense to 
do so for an aircraft target in view of the more dominant errors 
inherent in gunnery. One of my first reforms was simply to introduce 
4-figure trigonometric tables, and to equip the computing room with 
desk calculating machines in place of longhand pencil and paper sums. 
The calculating machines were winkled out of the Treasury, who were 
keeping them massed in a big cupboard in case they might be of 
future service for financial purposes.

There were certain bits of mathematics, of which I had no previous 
knowledge; in particular I needed to learn about numerical methods 
and statistics. I taught myself from Whittaker and Robinson's book 
about subjects such as finite differences and interpolation. 
To describe the trajectory of a shell, given the angle of elevation 
of the gun firing it, range tables of the sum were available in 
terms of the Cartesian coordinates of the shell at successive widely 
spaced intervals along its trajectory. It had not occurred to the 
compilers of the range tables that it would be more natural to 
represent this data in terms of polar coordinates; and, even when 
this was done there remained the non-trivial task of 
two-dimensional interpolation of this data. There is a result, 
due to Kolmogorov, that a continuous function of $d$ independent 
variables can be expressed in terms of a polynomial in $2d+1$ functions 
each of a single variable; but I did not know of this result until 
well after the war was over. Nevertheless I discovered for myself 
shortly after arriving at Lydstep that this result was explicitly 
true in the particular case $d = 2$ at least for the polar coordinate 
versions of $3.7"$ and $4.5"$ anti-aircraft guns. Accordingly we
recalculated the range tables of these guns in terms of quadrant 
elevations and tangent elevations; and were then able to complete 
the predicted trajectory using $1$-dimensional interpolations.

Acquaintance with statistical techniques was the other main gap 
in my previous mathematical education; and to cover this I 
obtained leave of absence to return to Cambridge for a few weeks. 
The first volume of M.\ G.\ Kendall's book on mathematical statistics 
had just been published. I also read R.\ A.\ Fisher's book on statistical 
methods for research workers. Statistical techniques played an 
important role at Lydstep in assuring the performance of anti-aircraft 
radars and predictors, and in liaising with radar developments 
from the Radar Research Establishment at Malvern.

By the end of the war I had been promoted to the rank of major, 
and appointed a consultant to the Ordnance Board in London. 
Anti-aircraft gunfire, which had been pretty inaccurate at the 
beginning of the war, had gradually improved by the end of the war; 
in particular the V1 bomb was comparatively easy to shoot down because 
of the introduction of the proximity fuse in shells. Against this, 
the V2 bomb was a ballistic missile and so unassailable. In the 
near future hostilities with nuclear weapons would render discussions 
with the Ordnance Board about the air defence of London nugatory. 
Effectively, the chapter on anti-aircraft gunnery was closed.
\endquote

\msection{Postwar Activities}
\quote
In 1946 I returned to Cambridge as an undergraduate at Emmanuel 
College. From time to time there were occasional trips up to London 
to fulfill my duties at the Ordnance Board, but these had 
little relevance to the future of anti-aircraft gunnery. 
Before the war I had done a certain amount of skiing; and I 
hoped for a half blue for skiing. One of the difficulties was 
that foreign currency was rationed by the Treasury; and so I 
needed to earn some Swiss francs by giving some lectures on 
statistics at any Swiss university that could be persuaded to 
employ me. Thanks to references provided by Harold Jeffreys, 
the Federal Institute of Technology (E.T.H.) in Zurich was kind 
enough to provide the necessary funds. However in those days the 
university skiing team consisted of four members, and I was ranked 
fifth in the trials; so I never got a half blue, although I did 
take part in a joint Oxford--Cambridge match against the combined 
Swiss universities which was a twelve-a-side match. 
Needless to say, the combined Swiss universities beat the joint  
Oxford--Cambridge team. 

As a Cambridge undergraduate in the two years after the war I was 
much more motivated than I had been in 1939/40; and I also had 
the good fortune to be tutored by better tutors, in particular 
A.\ J.\ Ward and J.\ A.\ Todd for pure mathematics and R.\ A.\ Lyttleton 
for applied mathematics. In 1948 I got a first class (Wrangler) 
in Part 2 of the Mathematical Tripos.

In 1948 I thought I would like to try my hand at an academic job 
in mathematics or mathematical statistics. There was no opening 
for me at Cambridge then. I applied for vacant lectureships at 
Reading University and at St Andrews University, but my applications 
were not successful. However I did get an appointment as a graduate 
assistant at Oxford in the Lectureship in the Design and 
Analysis of Scientific Experiment.

This Lectureship was a small department headed by the lecturer 
(D.\ J.\ Finney) and having two graduate assistants (M.\ Sandford and myself) 
together with a secretary and a couple of girls with desk calculators. 
At that time it was the only established provider of statistical 
services at Oxford, and its remit was spread quite generally over 
any and all queries that might be thrown up in various branches of 
service. It also had to offer lectures and instructions on statistics; 
for example, it fell to me to give the lecture course in the 
Department of Forestry for overseas forest officers on the 
collection and analysis of data on trees and their growth. 
\endquote

\msection{Oxford}
Hammersley held the position of Graduate Assistant, Design and Analysis of Experiments, at 
Oxford University, until he moved in 1955 to AERE Harwell as Principal Scientific Officer.
He returned to Oxford in 1959 as Senior Research Officer at the Institute of Economics 
and Statistics.  This was a position of roughly the same level as a University Lecturer 
but with neither formal teaching duties nor a linked College Fellowship/Tutorship.

It was during this period that he began an association with Trinity College which was to last for 
the rest of his life. When P.\ A.\ P.\ Moran left Oxford for the Australian National University at 
the end of 1951, Hammersley took over his tutorial duties at Trinity
as Lecturer in Mathematics.  
It was not until his election to a Senior Research 
Fellowship in 1961 that he became a fellow
of the college.  In 1969 he was promoted to (University) Reader in Mathematical Statistics, 
and was elected to a Professorial Fellowship at Trinity,  
two positions that he retained until his retirement in 1987. 
It is sometimes said that Hammersley was only
the second mathematics fellow at Trinity since its foundation in 1555, 
following in the footsteps of Thomas Allen (elected in 1564). 
He was in fact arguably the first such fellow.
In the late 16th century
all Trinity fellows were required to take the oath of
supremacy, an obligation that
Allen avoided by departing the College in 1571.
It was during that period and later
that Allen's mathematical activities developed, although, 
unlike Hammersley, he
is said to have written ``little and published nothing'' (see \cite{dnb}).

Despite the fact that Hammersley held no official position 
at the University between 
1955 and 1959, he took on his first four Oxford DPhil students in October 
1956.  He retained an office in Keble Road, and he seems to have spent a 
lot of his time there.
From 1959 until his retirement in 1987, 
he worked in what appeared to be splendid isolation in 
his office in the Institute of Economics and Statistics in St Cross Road.  As far as one could 
judge, apart from seeing graduate students and teaching a few Trinity undergraduates, he had 
his time free for research.

It was over Sunday lunch in Oxford shortly after his arrival
that he met Gwen Bakewell, who became his wife in 1951.  
Their first home in Longwall Street was soon replaced by Willow Cottage 
on the Eynsham Road, 
where their sons Julian and Hugo were born.

Although his university position was not in mathematics, he was a member of the subfaculty,
and he lectured and examined under its auspices. He gained a certain
notoriety for his high expectations of undergraduates. For example, one year
he offered a non-examinable lecture course on \lq Solving Problems' 
in which few students lasted very long. As a Finals 
Examiner in 1966, he set (or was at least blamed for) what was the most difficult set of compulsory 
papers in living memory.  1966 became known as the `year of the carrot' 
in honour of one  
question on differential equations that opened with the phrase: 
``If a sliced carrot is immersed at 
time $t = 0$ in $\beta$-indolyl acetic acid $\dots$''

Basic mathematical techniques mattered a lot more to Hammersley than many an
advanced result. On one occasion in an examiners' meeting, 
he would not withdraw from the
position that a relatively large number of marks, in an advanced probability
question, be given for the correct use of partial fractions.

It was not always easy for students and colleagues to rise to the uncompromisingly high 
intellectual standards set by John Hammersley, but it was a level playing field, and he 
applied his standards to himself just as to others.  To the knowledge of the current authors,
he took on only eight doctoral 
students during his career, and at least five of these continued 
to successful scientific careers.  Students were required to show their worth, as explained by John Halton:

\quote
A cousin drew my attention to an advertisement in the Observer
$\dots$, seeking applicants for UK Atomic
Energy Authority Research Studentships, to study Monte Carlo
methods for a DPhil at Oxford. $\dots$
In a few weeks, I was invited to ``present myself for examination"
at the UKAEA site at Didcot. With very little idea of what this
would entail, I went.
There I found a [number of] equally bemused applicants,
who were ushered into a large hall furnished with a suitable number
of small desks and sat down. John Hammersley strode breezily up
to the podium, introduced himself, and asked us to write a four-hour
examination, consisting of a dozen or so tough mathematical questions.
I attempted to solve each problem in turn, suggested possible lines
of approach, and tried to answer the questions posed, with little
success.
At the end of four hours, the papers were collected and we waited
anxiously for the outcome. 
\endquote

\comment
Much later, John Hammersley told us
that the questions asked on the four-hour examination had been research
questions that he had worked on but not solved, at that time!
I was emotionally crushed, sure that I had failed dismally. Later,
just nine of us were asked to come back before a board of scientists
for an interview. I was so shaken by then that, when asked what kinds
of problems I was working on at EE [English Electric], 
I simply could not recall!
Some days later, we were told that five of us \hbox{\rm[one later withdrew]}
had been selected
to receive the Research Studentships and to begin work
in the Michaelmas Term of 1956. I believe that we were
the first graduate students that John taught at Oxford $\dots$
\endcomment

Peter Marcer has taken up the story:
\quote
What a sleepless night I (and I expect others) had before the interviews the
next day, when each of us asked members of the panel, which included
John and Professor Flowers as he was then, what the answers were and how
one did the questions. Only to be told that John had done the rounds of
the theoretical physics department at Harwell, and compiled the
examination out of the questions that the members of that department were in
the course of trying to answer!  That is, there were no answers to these
questions as yet, and the panel just wanted to see how we, the candidates,
might begin to tackle them!  I think that episode sums up John for me, a
great mind sometimes puckishly inclined but with great purpose, and above
all a great gentleman of the old school. He was a delight to know, and
will be sorely missed, and I owe him a great deal.
\endquote

As a result of this exercise,
Halton, Marcer, David Handscomb, and Jillian Beardwood were awarded
studentships under Hammersley's supervision. As `Monte Carlo' students, they were
privileged with access to the Ferranti Mercury computers
at both Oxford and Harwell, as well as to the Illiac II while 
visiting the University of Illinois at Urbana in 1958.

Hammersley was for a period equally at home in California and Oxford.  
He was a regular 
contributor to the Berkeley Symposia on Mathematical Statistics and Probability, and was a 
close friend of the distinguished statistician Jerzy Neyman.  He spent the Michaelmas 
terms of 1958 and 1961 at Urbana, Illinois and Bell Telephone Laboratories, Murray Hill 
respectively.  On both these trips he was accompanied by his graduate students.

He never studied for a PhD, perhaps because 
of his age following war service, 
but he was awarded an ScD by Cambridge University in 1959,
followed in the same year by an Oxford DSc (by incorporation).  
He was awarded the Von Neumann 
Medal for Applied Mathematics by the University of Brussels (1966), 
the Gold Medal of the Institute of Mathematics and its 
Applications (1984), and the P\'olya Prize 
of the London Mathematical Society (1997).
He was elected to the Royal Society in 1976. 
He gave the 1980 Rouse Ball lecture at Cambridge University, and
published an account in \cite{82a}.

On retiring from his Oxford Readership in 1987, he was welcomed at 
the Oxford Centre 
for Industrial and Applied Mathematics (OCIAM).  
He reciprocated
this act of hospitality by making his extensive mathematical experience 
available to all who asked.

Many of Hammersley's friends and colleagues gathered in 1990 at the Oxford Mathematical 
Institute for a conference to recognise his 70th birthday.  A volume \cite{GW} entitled 
\lq Disorder in Physical Systems' was published in his honour, 
with contributions from many whose work 
had been touched by his ideas.  Hammersley gave the closing lecture of the meeting under 
the title \lq Is algebra rubbish?', but he uncharacteristically refrained on this occasion 
from answering the question.

In more recent years he was to be found at Willow Cottage, reading, doing the crossword, 
and working on Eden clusters.  He died on 2 May 2004 after an illness.

\msection{John Hammersley, mathematician}
John \h\ was an exceptionally inventive mathematician and
a remarkable and fearless problem-solver. 
He had the rare ability to pinpoint the basic mathematics
underlying a scientific problem, and to develop a useful theory.
He preferred what he called ``implicated" mathematics over
``contemplative" mathematics; that is, he found the solution
of problems to be superior to the ``high-rise mathematics"
of which he could be sharply critical (see [\citem{new1}, \citem{new3}]). 

The conventional modern classification of mathematics into pure, applied,
and statistics can accentuate gaps between these areas, gaps that
need to be filled.  Hammersley spurned such an attitude; when
facing a practical problem, he used whatever he could find in order to solve it.
This `bare hands' approach does not always lead to the neatest
solution, although, in Hammersley's case, much of the resulting mathematics
has stood the test of time. Several of the problems that he 
formulated and partly solved have emerged since as landmark problems
of combinatorics and probability. For example, his work on self-avoiding
walks and percolation is fundamental to the theory
of stochastic L\"owner evolutions that is now causing a re-think of
the relationship between probability and conformal field theory;
his results on the Ulam problem underly the proof \cite{BDJ}
that the relevant weak limit is the Tracy--Widom distribution. 
These two general areas are
amongst the liveliest of contemporary mathematics, as witnessed
by the award of Fields Medals in 2006 to Werner and Okounkov.

Paper \cite{22}, written jointly with K.\ W.\ (Bill) Morton, is a landmark
of his earlier work in two regards. Firstly, it marks a beginning
of Hammersley's extensive study of discrete problems in probability
and statistical mechanics. Secondly,
the paper contains two
problems and a technique that have attracted a great deal
of attention in the 50 years since.  Despite the title of the paper,
`Poor man's Monte Carlo', the lasting contributions are the
clear statement of the problem of counting self-avoiding
walks, the use of subadditivity to prove the existence of the
connective constant, and the discussion of
random media that culminated
in Simon Broadbent's formulation of the percolation model.

These and other topics are discussed further in the following paragraphs,
complemented by summaries of how John's work has stimulated the
relevant fields since. 

\ssection{Computing/Calculating/Estimating}
Hammersley's early scientific work was based on the mathematics he had been 
doing during the war.  His first publication \cite{1} arose from 
independent contributions by Majors Bayley and 
Hammersley to the discussion following the reading of a paper on random processes 
by Maurice Bartlett \cite{Bart} at a symposium 
on Autocorrelation in Time Series held in 1946 at the Royal 
Statistical Society.  The problem confronting Bayley and Hammersley arose 
in trials of anti-aircraft equipment.  
The details were embodied in ``reports not 
generally available'' but \cite{1} contains in condensed form some of the results obtained.

There followed a sequence of papers on essentially unrelated problems, many concerned with hard 
calculations or estimation. Probably his first significant work was his paper \cite{6} on the 
estimation of parameters when the parameter space is a discrete set of points.	
He showed, for example, that, 
if the unknown mean of a normal population with given variance is assumed integer-valued,
then its maximum likelihood estimator is the integer 
nearest to the sample mean.
His interest in issues of this kind arose from a problem of 
estimating the molecular weight of insulin, and this may 
have come to his attention during his work as a 
consultant on statistical problems to members of the university in the natural 
sciences.

It was a mathematical problem arising in \cite{6} that led to his paper 
\cite{9} on asymptotic formulae for the sums of products of the natural numbers.  
Paper \cite{9}, read in isolation, may appear to be scantily 
motivated.
However, it does display Hammersley's formidable analytic skills, and it 
attracted the attention of Paul Erd\H os who, in \cite{Erd}, settled one of the open 
problems posed.  It is now clear that, in \cite{9}, he
was in fact calculating what Cramer \cite{Cram} described recently as ``remarkable 
expressions'' for the mode of Stirling numbers of the first kind.

Throughout the rest of his scientific career, John Hammersley continued this 
interest in computing methods and computer science---principally through 
his work on large-scale simulations (see below).  

\ssection{Applied Probability}
In the period between leaving the military and starting his collaboration with 
Morton, Hammersley appears to have tried his luck at a range of problems in 
applied probability, hard analysis and large-scale computations.  For example,
in \cite{4}, he considers a problem arising in the design of experiments that
may be expressed as follows: given a collection of $k$ counterfeit and $n-k$ genuine coins,
how may we detect the counterfeit coins? 
His interest in stochastic geometry was developed in \cite{5}, 
where he studied the distribution of the distance between two points independently 
and uniformly distributed over the solid $n$-sphere. 
In \cite{10}, he proved a special case of a conjecture of Fejes 
T\'oth about the sum of the side-lengths of a convex polyhedron containing a sphere of unit 
diameter.  His paper \cite{20} on Markovian walks on crystal lattices originated from a 
study of diffusion of electrons in crystals such as the hexagonal close-packed lattice.

Around 1953, he considered a problem on counting blood cells 
that had arisen at the Clinical Pathology Department of the 
Radcliffe Infirmary at Oxford.  The mathematical problem here turns out to be 
equivalent to finding the probability distribution of the number of gaps between 
intervals of random length placed randomly on a circle.  Hammersley showed (by 
typically hard analysis) that it was asymptotically normal.  Cyril
Domb has given an account in \cite{Domb} 
of the history and ramifications of this particular problem,
and this work illustrates Hammersley's gift for picking out hard, 
genuinely interesting problems from the applied sciences and translating them 
into valid mathematics.

In \cite{30}, he extended a classical result of Mark Kac \cite{Kac}
concerning the number of zeros of a polynomial with random coefficients.
Kac's results were for the mean number of real zeros when the coefficients are 
independent, identically, and normally distributed, and Hammersley gave a 
substantial, albeit complicated, generalisation.  For recent activity
in this field, see [\citem{Farah}, \citem{Fried}, \citem{Ramp}].

Hammersley's most influential work in applied probability is that on percolation and
on the large-scale geometrical properties of
$n$ points dropped at random into a bounded region of Euclidean space.
We return to these two areas in later paragraphs. 

Having sketched Hammersley's early work, we move to his work post-retirement, almost
all of which was concerned with the growth of crystals.
He worked with Mazzarino on a 
third-order differential equation arising as a model for the growth 
of a crystal in a supercooled liquid, see [\citem{83}, \citem{84}].
This `classical' work was followed by his final two research papers
directed at the stochastic model introduced in 1961 by Murray Eden \cite{Eden} for 
growth in biological cells.  Despite its apparent simplicity, the Eden model has attracted 
a great deal of interest over the years.  

In the simplest version, the `cells' are taken to be closed 
unit squares of the 2-dimensional square lattice.  
All cells but one
are coloured white initially, and subsequently cells are blackened 
one at a time.  The mechanism of growth is as follows. An edge of the 
lattice is called {\it active\/} if it separates a black cell 
from a white cell.  At stage $n$, an active edge is picked at 
random, and the associated white cell is coloured 
black.  At time $n$, there is a cluster $C_n$ containing 
$n + 1$ black cells.  The 
shapes of the $C_n$ have the same distribution as those
of the first-passage percolation model
discussed below, when the edge-passage-times of that model
are exponentially distributed. 

Natural questions of interest about this process are: (i) what is 
the `shape' of $C_n$ for large $n$, and 
(ii) how large do the `lakes' of enclosed white cells grow before 
they are eventually filled in by black cells and disappear?
In \cite{90}, Hammersley presented non-rigorous arguments suggesting that 
all lakes in the `island' $C_n$ lie with high probability
within a distance $\O(\log n)$ of the coastline.
      
In his penultimate research paper \cite{89}, Hammersley (with Mazzarino)
carried out a large-scale Monte Carlo simulation in which clusters of 
size of order $10^9$ are grown, and various quantities such as the mean 
cluster-radius are estimated.  The authors evinced
pride in being able to carry out this huge computational task using only 24 megabytes 
of a Convex 220 machine, in contrast to comparable simulations 
of Zabolitsky and Stauffer \cite{ZabS} using a Cray 2 with 4 parallel processors
and a vast (for the period) store of 2045 megabytes.
      
A subject of primary interest in these two papers is the `surface roughness' of a typical cluster.  
The theoretical analysis carried out in \cite{90} makes use
of the theory of {\it harnesses\/}, as introduced by Hammersley in 1967.  
Harnesses may be described loosely as a spatial generalisation of a martingale;
they appear to have received very little attention since 1967, although 
Hammersley's original paper \cite{59} was 
one of the 45 articles selected and reprinted in \cite{FamV} 
as one of the seminal contributions to the scaling laws 
that characterise rough surfaces generated stochastically.

\ssection{Monte Carlo Methods}
From the very beginning of his career, John Hammersley sought
methods to carry out large computations. The equipment then
available was limited and unreliable and, rather as in his army days, he 
became a master of desk calculators and early computers.
He considered it a virtue to use computing resources
in an economic and efficient manner, and this attitude remained with
him all his life. He once boasted of holding the 1961 world record 
for keeping a computer (at Bell Labs) working without breakdown for 39 hours.

Credit for the name and the first systematic development of Monte Carlo 
methods is usually
accorded to Fermi, Metropolis, von Neumann, and Ulam. 
This area
fascinated Hammersley. The idea is that one may estimate a quantity
through computations involving random numbers.
A principal objective is to reduce the degree of variation in the estimate,
thereby improving the accuracy of the result.

Hammersley's interest in Monte Carlo methods seems to have been sparked
by his attendance at a symposium in Berkeley in the early 1950s, and he
gave a Master's level lecture course on the subject on his return to Oxford.  
In the audience was Bill Morton who had 
just graduated (in 1952) from Oxford and held an appointment at AERE Harwell.  
It was around this 
time that Hammersley organised the workshop on Monte Carlo methods at Harwell
during which he met Simon Broadbent.

It was with Morton that Hammersley wrote his paper \cite{22} entitled `Poor man's Monte Carlo',  
of which the basic thesis was that one does not necessarily need large high-speed machines 
to use Monte Carlo effectively. 
In order to illustrate this main point, the authors draw on a range of examples such 
as self-avoiding walks.  Among the more diverting of the 
examples is the testing of a quantum hypothesis of Alexander Thom.
Thom had measured the diameters of 33 Druid circles in Western Scotland, and, based on 
the (integer) data, he conjectured that these diameters were intended to be multiples 
of 11.1 feet.  The evidence for this was that 27 of the circles had diameters lying in the range 
$11.1(n \pm \frac14)$ for integral $n$.
Hammersley and Morton used simple Monte Carlo methods to test the 
hypothesis and, as David Kendall suggested in \cite{DGK}, their work led to a statistical examination 
which went a long way towards confirming this proposal.

Monte Carlo methods are based upon the use of pseudo/quasi-random
numbers, and this raises certain issues of principle. 
Hammersley's impatience with philosophical discussions involving the ethics or 
correctness of using pseudo/quasi-random numbers in place of truly random ones  
is captured in his reply to the discussions at the Symposium on Monte Carlo 
Methods at which \cite{22} was presented: 
``The discussion has raised several questions about random numbers: do they even exist; 
can they be produced to order and if so how; can they be recognised and can we test that 
they are not imposters?  These are diverting philosophic speculations; but the applied 
mathematician must regard them as beside the point.''

Indeed, his intolerance of philosophy as an academic subject seemed to stay with him 
throughout his life.  The Oxford
joint school of Mathematics \& Philosophy was one of his 
{\it b\^etes noires\/}, and various amusing stories have accumulated about the year where he  
ended up (by default) as Chairman of the Examiners.  When the
opportunity came for him to chair the Finals examining board, he grasped it
enthusiastically, and taped (with his colleagues' permission) a post-meeting discussion on
the value of the degree. His further strenuous
efforts could not in the end persuade either the mathematicians or the
philosophers that the degree should be shelved.

Hammersley's most significant contribution to the theory, as against practice, of 
Monte Carlo methods is probably his work on antithetic variates.  This is a technique 
for yielding estimates with variances considerably less than those obtainable by a 
naive approach.  This is typically achieved by representing the estimator 
as a sum of correlated random variables, and it is one of the most popular variance-reduction 
techniques.  Its drawback is that many 
antithetic sampling plans are too computationally complex to be of practical use 
in simulations.  Despite this, the work of Hammersley and Morton \cite{28} is currently 
regarded as a major contribution.  
(See, for example, \cite{RoaW}.)  
It is interesting, therefore, that in \cite{53} 
Hammersley and Handscomb claim only the name, not the original idea which, as 
pointed out by Tukey \cite{Tukey}, can be regarded as an important special case of regression.
This technique is now, perhaps, one of the most important in the 
application of Monte Carlo methods to high-dimensional numerical integration,
with applications in many areas including mathematical finance.

The Hammersley--Handscomb monograph \cite{53}, published in 1964, is a 
landmark in the study of Monte Carlo methods and is still much used today.  
Hammersley's interest in the field seems to have declined 
following its publication.  

\ssection{Percolation}
Percolation was born as a mathematical object
out of the musings on random media
found in \cite{22}, and it has emerged as a cornerstone of
stochastic geometry and statistical mechanics. 
One of the discussants of \cite{22}, Simon Broadbent,
worked at the British Coal Utilization Association, where
he was involved in the design of gas masks for coal miners
(see [\citem{opt}, \citem{78}]).
Hammersley recognised the
potential of Broadbent's proposal for flow through a random medium,
and they collaborated on 
the seminal paper \cite{31}, where the critical percolation
probability was defined.  There are earlier references 
to processes equivalent to percolation, see \cite{dVW} for example,
but it was Hammersley who initiated a coherent mathematical theory.

The basic model is as follows.  Consider a crystalline
lattice. We
declare each edge of the lattice (independently) 
to be {\it open\/} (to the passage of fluid) with
probability $p$, and otherwise {\it closed\/}. Fluid is
supplied at the origin of the lattice
and allowed to flow along the open edges only.
The fundamental question is to describe the size and geometry of the
set $C$ of vertices reached by the fluid.  
The significance of this  model is
far-reaching in stochastic geometry
and statistical mechanics, and the associated 
mathematics and physics literature
is now very extensive indeed. Of primary importance is the existence of
a phase transition: there exists a critical value $\pc$ such that
$C$ is finite when $p<\pc$, and $C$ is infinite with a strictly
positive probability when $p> \pc$. The non-triviality of the
phase transition was proved by Hammersley, as follows. 
Hammersley and Broadbent \cite{31}
established a lower bound for $\pc$ in terms of counts of
self-avoiding walks and the connective constant. (An account of
the connective constant may be found in the next section.)
This result was strengthened in \cite{33}, where it was shown that
$|C|$ has an exponentially
decaying tail whenever it has
finite expectation. The method developed in \cite{33}
is a precursor
of a now standard argument attributed to
Simon and Lieb [\citem{Li80}, \citem{Si80}] and usually expressed as: finite
susceptibility implies exponentially decaying correlations.
In \cite{34}, he proved an upper bound for $\pc$ in terms
of the boundary sizes of neighbourhoods of the origin,
and he deduced by graphical duality that $\pc<1$ for 
oriented and unoriented percolation on the square grid; this 
is the percolation equivalent of the Peierls argument for the Ising 
model, \cite{Peierls}. This general route to showing the 
existence of a phase 
transition is now standard for many models.

In an alternative model, it is the vertices 
rather than the edges of the crystal lattice
that are declared open/closed.  Hammersley \cite{40}
proved the useful fact that
$C$ tends to be smaller for the `site' model than for
the `bond' model, thereby extending a result of Michael Fisher. 
The best modern result of this type 
is by one of his students, see \cite{GrS96}.

An inveterate calculator, Hammersley wanted to calculate or
estimate the numerical value of $\pc$ for the square grid.
Theodore Harris proved in 
a remarkable paper \cite{Ha60} that $\pc\ge \frac12$, and Hammersley's
numerical estimates indicated $\pc<\frac12$; 
``what better evidence could exist for $\pc=\frac12$",
he would ask. He was therefore thrilled when Harry Kesten, \cite{K80},
proved the holy grail. This was however only the end of the beginning
for percolation.

\topinsert
\centerline{\epsfxsize=11cm \epsfbox{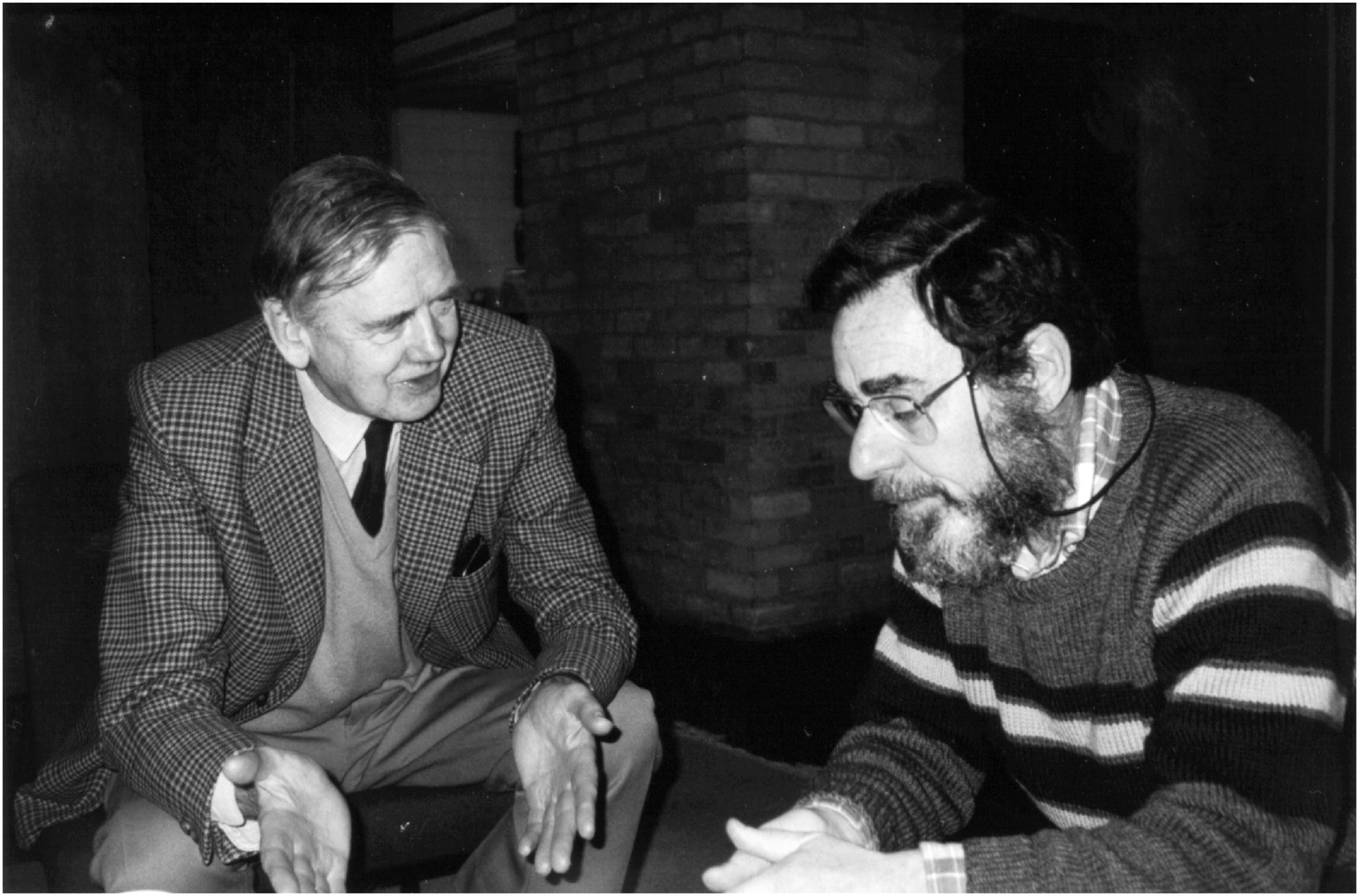}} 
%\centerline{\vbox{(photo here)}}  
\botcaption{John Hammersley and Harry Kesten in the 
Mathematical Institute, Oxford University, November 1993.}
\endcaption
\endinsert
   
Percolation Theory has gone from strength to strength in recent years.
The main questions are largely solved (see \cite{G99}), 
and current attention is focussed on the nature
of the phase transition in two dimensions.  Schramm
\cite{Sch00} predicted that the scaling limit of the perimeters of large
critical percolation clusters constitute a stochastic L\"owner evolution
(SLE)\footnote{Often termed a {\it Schramm--L\"owner evolution}.} with parameter 6.  
Smirnov \cite{Smi0} proved Cardy's formula
for crossing probabilities of critical site
percolation on the triangular lattice, and indicated how to
achieve the full scaling limit. See \cite{Sch06} for a survey of
SLE and associated problems and conjectures.

\ssection{Self-Avoiding Walks and the Monomer--Dimer Problem}
In the paradigm of statistical mechanics, a system is
modelled by a set of configurations to each of which is allocated a weight.
The sum of all weights is called the `partition function' and
the state of the system may be described via an analysis
of this function and its derivatives. In a system of
polymers, the first calculation is to find the number
of such polymers. When the polymers are simple chains rooted at the
origin of a lattice, this is the problem of counting self-avoiding walks
(SAWs). Let $s_n$ be the number of
SAWs of length $n$ on a given lattice.  The first serious progress towards
understanding the asymptotics of $s_n$ as $n\to\infty$ was made
in \cite{22}. The key is the `subadditive inequality'
$t_{m+n} \le t_m+t_n$ satisfied by $t_n=\log s_n$, 
from which the existence of the 
so-called connective constant
$\kappa=\lim_{n\to\infty} n^{-1} \log s_n$ follows immediately.
This observation, regarded now as essentially trivial in the light of
the complicated analysis achieved since,  has had a very substantial
impact on spatial combinatorics and probability.  It marked the
introduction of subadditivity as a standard tool, and it 
initiated a detailed study, still ongoing, of the geometry
of typical instances of geometrical configurations such as paths
and lattice animals.

The subadditive inequality implies the bound $s_n \ge \kappa^n$.
Hammersley put a great deal of energy into trying to find a 
complementary upper
bound on $s_n$,
but with only partial success. With his student Welsh,
he proved in \cite{46} that $s_n \le \kappa^n \exp(\lambda n^{1/2})$ for
some $\lambda<\infty$. This was improved by Kesten \cite{K7} for $d \ge 3$,
and such bounds were the best available for some time before
it was realised by others that a lace expansion could
be used for sufficiently high dimensions, see \cite{MadS}. 
As a result of a large amount of hard work and
some substantial mathematical machinery,
the problem of counting SAWs was solved by Hara
and Slade \cite{HS92} in five
and more dimensions. The case of two dimensions,
for which the bound of \cite{46} remains the best known,
has attracted much interest in recent years with the
introduction by Schramm of stochastic L\"owner evolutions (SLE),
and the conjecture that a random SAW in 
two dimensions converges in an appropriate sense as $n\to\infty$  to
a SLE with parameter $\frac 83$, see \cite{Sch06}. This conjecture
is one of the most important currently
open problems in probability.

Hammersley was happy in later life
to learn of progress with percolation and self-avoiding walks.
He felt that he had ``helped them into existence"
for others to solve. The two-dimensional percolation and SAW problems  
are two of the hottest
problems of contemporary probability, in testament to Hammersley's excellent scientific taste. 

There is a second counting problem of statistical mechanics that
attracted Hammersley, namely the {\it monomer--dimer\/} problem.
This classical 
problem in solid-state chemistry may be formulated   
as follows:
A {\it brick\/} is a $d$-dimensional ($d \geq 2$) rectangular parallelepiped 
with sides of integer lengths and even volume.  A unit
cube is called a {\it monomer\/}, 
and a brick with volume 2 a {\it dimer\/}. 
The dimer problem is to determine the number $f(a_1,a_2, \dots, a_d)$ of dimer 
tilings of the brick with sides of length $a_1,a_2, \dots, a_d$.  Hammersley 
proved in \cite{55} that the sequence $(a_1 a_2\cdots a_d)^{-1} \log f(a_1,a_2, \dots, a_d)$ 
approaches a finite limit $\lambda_d$ as the $a_i \rightarrow \infty$,
but what is the numerical value of $\lambda_d$?  
There is a `classical' result of statistical physics of Temperley--Fisher 
\cite{TempF} and Kasteleyn \cite{Kast},
who showed independently in 1961 that $\lambda_2$ exists
and is given by $\lambda_2 = \exp(2G/\pi) = 0.29156\dots$ 
where $G$ is Catalan's constant.
Hammersley devoted much energy on theoretical and computational approaches to  
finding a corresponding result for $d \ge 3$ but, as far as we know,
the exact value is still unknown even when $d=3$.  

In its more general form, the monomer--dimer problem amounts to the purely 
combinatorial question of counting the number $f_G(N_1,N_2)$ of distinct arrangements of 
$N_1$ monomers and $N_2$ dimers on the edges and vertices of a graph $G$, 
such that each dimer is placed on an edge, each monomer on a vertex, and 
each vertex of $G$ either is occupied by exactly one monomer or is the 
endvertex of exactly one dimer.  For this to be possible, $G$ must have 
exactly $N_1 + 2N_2$ vertices, and the {\it density\/} $p$ of the configuration 
is defined as the ratio $2N_2 / N_1$.
Hammersley proved in \cite{55} that
the number of $p$-density configurations on the cube of volume $n$
in $d$ dimensions is of order
$\lambda(d, p)^n$ for some function $\lambda$.
He spent much effort on obtaining bounds for $\lambda$, but, 
even today in 2-dimensions, our knowledge
is very limited.  See, for example, \cite{FriedP}.

The dimer problem is very much alive today. The
two-dimensional model turns out
to be related to the Gaussian free field and to stochastic L\"owner evolutions
with parameters 2, 4, and 8,
see [\citem{Ken1}, \citem{KenW}] for example.

\ssection{First-Passage Percolation, and Subadditive Processes}
Percolation is a static model in the sense that each edge
is either open or closed, and water is considered to flow 
instantaneously along open edges. Hammersley and Welsh formulated
a time-dependent version of this model in \cite{52}, and dubbed
this `first-passage percolation'. To each edge of the lattice is assigned
a random passage-time, and the time $a_{x,y}$ for water to
reach a given point $y$, having started at 
$x$,
is the infimum over all 
paths $\pi$ from $x$ to $y$ of the aggregate passage-time
of edges in $\pi$.  This pioneering paper \cite{52}
is now recognised
as one of the first works of mathematical
significance in the theory of the spread of material,
whether it be disease, fluid, or rumour,
through a random medium. The basic problem was to prove
the existence of a speed function $\sigma_{x} = \lim_{n\to\infty}
n^{-1} a_{0,nx}$, where 0 denotes the origin of the lattice. 
Hammersley and Welsh realised
that the key lay in the use of subadditivity,  $a_{0,x+y} \le a_{0,x}
+a_{x,x+y}$, the difference with previous applications being
that this inequality involves {\it random variables\/} rather than 
deterministic quantities.

They proved a version of the subadditive limit theorem
for stationary stochastic processes indexed by $d$-dimensional space,
the first `subadditive ergodic theorem'. They realised that this is best
done in the context of a general set of assumptions, rather than the specific
situation outlined above, and thus their paper gave birth to
one of the principal techniques for the analysis of
spatial random processes.  The search began for the `right'
combination of definition/theorem, and this was found by 
John Kingman \cite{K1} in one of the classic papers of 20th
century probability.  Despite later elaborations, it remains
fascinating to read this early literature, and especially
the dialogue of [\citem{K2}, \citem{65}].
Kingman's invited review article \cite{K2} (with
published discussion) appeared in the Annals of Probability. 
Hammersley's  
contribution to this discussion was too extensive to be accepted
as such by the journal editor, and it appeared later as \cite{65}.
It is there
that the condition of pathwise subadditivity is
replaced by the weaker assumption of `superconvolutivity'
of the associated probability measures.

In an earlier application of subadditivity to spatial systems
pursued jointly with his students Jillian Beardwood and John Halton,
Hammersley made a fundamental contribution to
the study of typical instances of problems in
operations research. Drop $n$ points at random into
a plane region $R$ of finite area. What is the length of the minimal
spanning (Steiner) tree and the minimal 
travelling salesman path on these points?
They showed in their classic
paper \cite{35} that the answer is (in essence) proportional to
$c_R\sqrt {n}$ for some constant $c_R$, and they developed also a higher-dimensional theory.
The key was to encode the problem in such a way that the
natural length-scale is $\sqrt n$, and then to use
a type of spatial subadditivity. This theorem was
central to the later work of Karp on a probabilistic
analysis of the random Euclidean travelling salesman problem, \cite{K77}. 
Further developments are described
in the Festschrift paper \cite{S90} by Steele. 

The title of \cite{S90} makes play on Hammersley's own
famous title `A few seedlings of research', published in 1972 in
the Proceedings of the Sixth Berkeley Symposium. In this
inspiring account of how to do mathematical research, Hammersley
showed in particular how to use subadditivity to solve (in part) the
now famous Ulam problem: in a random permutation of the first
$n$ natural numbers, what is the length $l_n$ of the longest
increasing subsequence? It turned out for geometrical reasons
related to \cite{35} that the answer is 
asymptotically $c\sqrt n$. This was the starting point of 
a major area of probability theory. Hammersley
claimed a back-of-the-envelope argument to show $c=2$,
but the formal proof eluded him and was found by
Vershik and Kerov \cite{VerK} and Logan and 
Shepp \cite{LogS} in the context of
random Young tableaux.  Interest turned then to
the size of the deviation $l_n-2\sqrt n$. Many partial
results preceded the remarkable proof by 
Baik, Deift, and Johansson \cite{BDJ}
that $(l_n-2\sqrt n)n^{-1/6}$ converges
as $n\to\infty$ to
the famous Tracy--Widom distribution of random matrix
theory.

\ssection{Random Fields}
One of the most important topics in modern statistics is the 
Bayesian theory of image
analysis. In this study of spatial random systems,
it is useful to have a classification of those probability
measures that satisfy a certain 
`spatial Markov property', namely: the configuration inside any region
$V$ depends on the configuration outside $V$ only through
the states of the vertices on its boundary $\partial V$. Some limited theory
of such measures was developed by Averintsev,
Dobrushin, Spitzer, and others around 1970.
This was generalised to an arbitrary network by Hammersley in 1971
following a suggestion of Clifford (see [\citem{Cliff}, \citem{HamC}]). 
The ensuing 
theorem,
commonly termed the Hammersley--Clifford theorem, though never
formally published, is much used in probability and statistics. It states that
a positive measure is a Markov field if and only if it
has a Gibbsian representation in terms of some potential
function. The methods used by Hammersley were much clarified by later 
authors including another of his students, Grimmett \cite{GMark},
who reduced the proof to an exercise in the inclusion--exclusion principle.

In Michaelmas
Term 1971, Hammersley offered a graduate course on Markov fields at the Mathematical Institute.
He promised a solution to the corresponding problem in which
the assumption of positivity is relaxed. It was typical of the man that
he had not proved the result, and indeed the `theorem' was
disproved through the discovery of a counterexample by a Rhodes Scholar,
John Moussouris, in the audience.
See [\citem{Mouss}, \citem{new3}].

\ssection{Educational Issues}
Great changes were made during John Hammersley's lifetime in the teaching of 
mathematics in schools, and he was for a period at the
forefront of the debate.  From the 1950s onwards, he argued fiercely
that schoolchildren and undergraduates
should be trained to solve problems, and that the curriculum should
be designed accordingly. He lectured 
on this topic around the UK, and he contributed 
to the development of the School Mathematics Project (SMP).
Not being a man of equivocal views, 
his uncompromising stance was seen by some as a provocation, but
he had many supporters and admirers.  However, the SMP proved no panacea 
for him: while it `modernised' aspects of mathematical teaching,
it introduced abstract theory without a sufficient problem element.

Hammersley frequently published his lectures in the Bulletin
of the Institute of Mathematics and its
Applications (IMA).  His principal article \cite{new1} on mathematical education 
appeared thus under the title 
`On the enfeeblement of mathematical skills by `Modern Mathematics' 
and by similar soft intellectual trash in schools and universities'.
This serious, if typically prolix, {\it critique\/} of school mathematics 
compelled a rebuttal from Bryan Thwaites, \cite{new54},
tempered as follows: ``I have, however, a profound reluctance to 
[reply to Hammersley's ``charges"]. 
The reason is that my admiration of the man and my opinion of his paper are in
great conflict. Much of my admiration stems from his mathematical
achievements; but it also rests firmly on my judgement
that it was he, more than any other Englishman, who
finally set going the long-overdue reforms in school mathematical
curricula."

Through his `popular' articles, Hammersley expressed his powerfully
held views on many matters, primarily scientific and educational.
These writings are erudite, provocative, and skilful with language,
if sometimes self-indulgent. His thoughts
on mathematical research were published alongside 
those of Michael Atiyah in [\citem{new55}, \citem{new3}], and include
some notable expressions: 
``$\dots$ perfuse his professorial piddledom", 
``Pure mathematics is subject to two diseases,
resulting from rigour and from axiomatisation", ``whatever
algebra can accomplish, some other branch of mathematics ought to be
able to accomplish more elegantly", ``$\dots$ and the production
of neater solutions is merely a matter for theory builders".
He loved a good phrase, even (perhaps, especially)
when it risked going a bit too far.
In reality, he would accept any theory that proved its worth.

As Hammersley wrote to Atiyah in \cite{new3}:
``I don't quarrel, but I am prepared to enter the lists. $\dots$
it is the jostling and jousting between different sorts of 
mathematicians and scientists, between different temperaments and 
unlike tastes,
that advances knowledge as a whole. So much the more fun, 
variety is the spice, and so on!" 

\msection{Acknowledgements}
We thank the Hammersley family for permission to
quote (in italics, and with minor corrections
and changes of presentation) from
John's account of his early life, written apparently in
response to a request from the Royal Society dated 1994. 
Christopher Prior and Clare Hopkins (archivist, Trinity College)
have advised us on College matters; John Halton
and Peter Marcer have reminisced about their experiences as
PhD students of John Hammersley. 
We thank Peter Collins, David Handscomb, and Bill Morton for their memories
of Hammersley in Oxford, 
and Harry Kesten for kindly commenting
on this biographical memoir.  

\def\PM{Philosophical Magazine}
\def\JSP{Journal of Statistical Physics}
\def\CMP{Communications in Mathematical Phy\-sics}

\def\TPR{The Physical Review}

\def\AP{Annals of Probability}
\def\JMP{Journal of Mathematical Physics}
\def\JAP{Journal of Applied Probability}
\def\AAP{Advances in Applied Probability}
\def\AAM{Advances in Applied Mathematics}
\def\JRSS{Journal of the Royal Statistical Society}
\def\JRSSB{Journal of the Royal Statistical Society B}
\def\JRSSA{Journal of the Royal Statistical Society A}
\def\PCPS{Proceedings of the Cambridge Philosophical Society}
\def\CM{Compositio Mathematica}
\def\AM{Acta Mathematica}
\def\JSIAM{Journal of the Society for Industrial and Applied Mathematics}
\def\BIMA{Bulletin of the Institute of Mathematics and
its Applications}
\def\AMM{The American Mathematical Monthly}
\def\MOR{Mathematics of Operations Research}
\def\SMD{Soviet Mathematics Doklady}
\def\JAMS{Journal of the American Mathematical Society}
\def\BLMS{Bulletin of the London Mathematical Society}
\def\IJM{Israel Journal of Mathematics}
\def\CRAS{Comptes Rendus des S\'eances de l'Acad\'emie des Sciences. S\'erie I.
Math\'ematique}
\def\DAM{Discrete Applied Mathematics}
\def\JLMS{Journal of the London Mathematical Society}
\def\BAMS{Bulletin of the American Mathematical Society}
\def\NAWIM{Indagationes Mathematicae}
\def\JCTB{Journal of Combinatorial Theory B}

\Refst\widestnumber\no{100}
\reset

\subrefhead{Book}

\ref
\noOf{53}
\by \jmh, Handscomb, D.\ C.
\book
Monte Carlo Methods
\publ Methuen
\publaddr London
\yr 1964
\endref

\subrefhead{Scientific articles}

\ref
\noOf{1}
\by Bayley, G.\ V., \jmh
\paper The effective number of independent observations in an 
autocorrelated time series  
\jour \JRSS\ (Supplement)
\vol 8  
\yr 1946 
\pages 184--197
\endref

\ref
\noOf{n1}
\by\jmh
\paper A geometrical illustration of a principle of
experimental directives
\jour \PM
\vol 39
\yr 1948
\pages 460--466
\endref

\ref
\noOf{n2}
\by\jmh
\paper An elementary introduction to some inspection procedures
\yr 1948
\vol 17
\pages 315--322
\jour Rev.\ Suisse Org.\ Indust.
\endref

\ref
\noOf{2}
\by  \jmh
\paper The unbiased estimate and standard error of the interclass variance  
\jour Metron  
\vol 15  
\yr 1949
\pages 189--205
\endref

\ref
\noOf{3}
\by \jmh
\paper
The numerical reduction of non-singular matrix pencils 
\jour Philosophical Magazine 
\vol 40 
\yr 1949
\pages 783--807
\endref

\ref
\noOf{n3}
\by\jmh
\paper Electronic computers and the analysis of stochastic processes
\jour Mathematical Computing
\yr 1950
\vol 4
\pages 56--57
\endref

\ref
\noOf{n4}
\by\jmh
\paper Calculating machines
\jour Chambers' Encyclopedia
\yr 1950
\endref

\ref
\noOf{n5}
\by\jmh
\paper Harmonic analysis
\jour Chambers' Encyclopedia
\yr 1950
\endref

\ref
\noOf{4}
\by \jmh
\paper
Further results for the counterfeit coin problems
\jour \PCPS
\vol 46
\yr 1950
\pages 226--230
\endref

\ref
\noOf{5}
\by \jmh
\paper
The distribution of distance in a hypersphere
\jour Annals of Mathematical Statistics
\vol 21 
\yr 1950
\pages 447--452
\endref

\ref
\noOf{6}
\by \jmh
\paper On estimating restricted parameters
\jour \JRSSB 
\vol 12
\yr 1950
\pages 192--229
\endref

\ref
\noOf{7}
\by \jmh
\paper
A theorem on multiple integrals
\jour \PCPS
\vol 47
\yr 1951
\pages 274--278
\endref

\ref
\noOf{8}
\by \jmh
\paper
On a certain type of integral associated with circular cylinders
\jour Proceedings of the Royal Society, Series A
\vol  210
\yr 1951
\pages 98--110
\endref

\ref
\noOf{9}
\by \jmh
\paper
The sums of products of the natural numbers
\jour Proceedings of the London Mathematical Society
\vol 1
\yr 1951
\pages 435--452
\endref

\ref
\noOf{10}
\by \jmh
\paper
The total length of the edges of the polyhedron
\jour Compositio Mathematica
\vol 9
\yr 1951
\pages 239--240
\endref

\ref
\noOf{n6}
\paper The absorption of radioactive radiation in rods
\by\jmh
\jour National Bureau of Standards W.P.
\vol 1929
\yr 1951
\pages 1--11
\endref

\ref
\noOf{11}
\by \jmh
\paper
The computation of sums of squares and products on a desk calculator
\jour Biometrics
\vol 8
\yr 1952
\pages 156--168
\endref

\ref
\noOf{11a}
\by \jmh
\paper
An extension of the Slutzky--Fr\'echet theorem
\jour Acta Mathematica
\vol 87
\yr 1952
\pages 243--257
\endref

\ref
\noOf{12}
\by \jmh
\paper
Lagrangian integration coefficients for distance functions taken over right
circular cylinders
\jour Journal of Mathematical Physics
\vol 31
\yr 1952
\pages 139--150
\endref

\ref
\noOf{13}
\by \jmh
\paper 
Tauberian theory for the asymptotic forms of statistical frequency functions
\jour \PCPS
\vol 48
\yr 1952
\pages 592--599
\moreref Corrigenda: {\bf 48} (1953), 735
\endref

\ref
\noOf{14}
\by \jmh
\paper
On a conjecture of Nelder
\jour \CM
\vol 10
\yr 1952
\pages 241--244
\endref

\ref
\noOf{n8}
\by\jmh
\paper Capture--recapture analysis
\jour Biometrika
\yr 1953
\vol 40
\pages 265--278
\endref

\ref
\noOf{15}
\by \jmh
\paper
Tables of complete elliptic integrals
\jour Journal of Research of the National Bureau of Standards
\vol 50
\yr 1953
\pages 43
\endref

\ref
\noOf{16}
\by \jmh
\paper
A non-harmonic Fourier series
\jour \AM
\vol 89
\yr 1953
\pages 243--260
\endref

\nref
\noOf{17}
\by \jmh
\paper
Corrigenda to the paper, Tauberian theory for the asymptotic forms of
statistical frequency functions
\jour \PCPS
\vol 49
\yr 1953
\pages 735
\endref

\ref
\noOf{18}
\by \jmh
\paper
On counters with random dead time I
\jour \PCPS 
\vol 49
\yr 1953
\pages 623--637
\endref

\ref
\noOf{19}
\by Antosiewicz, H.\ A., \jmh
\paper
The convergence of numerical iteration
\jour The American Mathematical Monthly 
\vol 60
\yr 1953
\pages 604--607
\endref

\ref
\noOf{20}
\by \jmh
\paper
Markovian walks on crystals
\jour \CM 
\vol 11 
\yr 1953
\pages 171--186
\endref

\ref
\noOf{n9}
\paper The consistency of stop-watch time-study practitioners
\by\jmh
\jour Occupation Psychology
\vol 28
\yr 1954
\pages 61--76
\endref

\ref
\noOf{22}
\by \jmh , Morton, K.\ W.
\paper
Poor man's Monte Carlo
\jour \JRSSB 
\vol 16
\yr 1954
\pages 23--38
\endref

\ref
\noOf{21}
\by \jmh , Morton, K.\ W.
\paper
Transposed branching processes
\jour \JRSSB 
\vol 16
\yr 1954
\pages 76--79
\endref

\ref
\noOf{23}
\by \jmh, Morton, K.\ W.
\paper
The estimation of location and scale parameters from grouped data
\jour Biometrika 
\vol 41
\yr 1954
\pages 296--301
\endref

\ref\
\noOf{n10}
\by Eyeions, D.\ A., \jmh, Owen, B.\ G., Price, B.\ T., Wilson, J.\ G.,
Morton, K.\ W.
\paper The ionization loss of relativistic mu-mesons in neon
\jour Proceedings of the Physical Society (A)
\vol 68
\yr 1955
\pages 793--800
\endref

\ref
\noOf{24}
\by \jmh
\paper
Storage problems
\jour Mathematische Annalen 
\vol 128
\yr 1955
\pages 475--478
\endref

\ref
\noOf{25}
\by \jmh , Nelder, J.\ A. 
\paper
Sampling from an isotropic Gaussian process
\jour \PCPS
\vol 51 
\yr 1955
\pages 652--662
\endref

\ref
\noOf{26}
\by \jmh
\paper
The area enclosed by P\'olya's walk
\jour \PCPS  
\vol 52 
\yr 1956
\pages 78--87
\endref

\ref
\noOf{n11}
\by\jmh
\paper Percolation in crystals: gravity crystals
\jour UKAEA TP
\vol 13
\yr 1956
\endref

\ref
\noOf{27}
\by \jmh
\paper
Conditional Monte Carlo
\jour Journal of the Association for Computing Machinery 
\vol 3 
\yr 1956
\pages 73--76
\endref

\ref
\noOf{28}
\by \jmh , Morton, K.\ W. 
\paper
A new Monte Carlo technique: antithetic variates
\jour \PCPS 
\vol 52 
\yr 1956
\pages 449--475
\endref

\ref
\noOf{29}
\by \jmh , Mauldon, J.\ G.
\paper
General principles of antithetic variates
\jour \PCPS 
\vol 52
\yr 1956
\pages 476--481
\endref

\ref
\noOf{30}
\by \jmh
\paper
The zeros of a random polynomial
\jour Proceedings of the Third Berkeley Symposium
on Mathematical Statistics and Probability
\vol II
\yr 1956
\pages 89--111
\endref

\nref
\noOf{n13}
\book Proceedings of the Oxford Mathematical Conference
for Schoolteachers and Industrialists
\ed \jmh
\yr 1957
\publ Times Publishing Comapny
\publaddr London
\endref

\nref
\noOf{n14}
\by\jmh
\paper The value of mathematics and its teachers
\jour ibid
\yr 1957
\endref

\ref
\noOf{31}
\by Broadbent, S.\ R., \jmh 
\paper
Percolation processes. I. Crystals and mazes
\jour \PCPS
\vol 53 
\yr 1957
\pages 629--641 
\endref

\ref
\noOf{32}
\by \jmh
\paper
Percolation processes. II. The connective constant
\jour \PCPS
\vol 53 
\yr 1957
\pages 642--645
\endref

\ref
\noOf{33}
\by \jmh
\paper
Percolation processes: Lower bounds for the critical probability
\jour Annals of Mathematical Statistics
\vol 28
\yr 1957
\pages 790--795
\endref

\ref
\noOf{n15}
\by\jmh
\paper Discussion on renewal theory and its ramifications
\jour\JRSSB
\vol 20
\yr 1958
\pages 287--291
\endref

\ref
\noOf{n16}
\by Egelstaff, P.\ A., \jmh, Lane, A.\ M.
\paper Fluctuations in slow neutron average cross-sections
\jour Proceedings of the Physical Society (A)
\vol 71
\yr 1958
\pages 910--924
\endref

\ref
\noOf{34}
\by \jmh
\paper
Bornes sup\'erieures de la probabilit\'e critique dans un processus de 
filtration 
\inbook Le Calcul des Probabilit\'es et ses Applications
\yr 1959
\publ CNRS
\publaddr Paris
\pages 17--37
\endref

\ref
\noOf{35}
\by Beardwood, J., Halton, J.\ H., \jmh
\paper 
The shortest path through many points
\jour \PCPS 
\vol 55
\yr 1959
\pages 299--327
\endref

\ref
\noOf{36}
\by \jmh
\paper
Monte Carlo methods for solving multivariable problems
\jour Annals of the New York Academy of Sciences
\vol 86
\yr 1960
\pages 844--874 
\endref

\ref
\noOf{37}
\by \jmh
\paper
Limiting properties of numbers of self-avoiding walks
\jour Physics Review 
\vol 118
\yr 1960
\pages 656
\endref

\nref
\noOf{n17}
\by Coulson, C.\ A., \jmh
\paper The bottleneck in British science and technology
\jour New Scientist
\vol 10
\yr 1961
\pages 499--500
\endref

\ref
\noOf{38}
\by \jmh
\paper
The number of polygons on a lattice
\jour \PCPS 
\vol 57
\yr 1961
\pages 516--523
\endref

\ref
\noOf{39}
\by \jmh
\paper
A short proof of the Farahat--Mirsky refinement of Birkhoff's theorem on
doubly-stochastic matrices
\jour \PCPS 
\vol 57
\yr 1961
\pages 681
\endref

\ref
\noOf{40}
\by \jmh
\paper
Comparison of atom and bond percolation processes
\jour Journal of Mathematical Physics 
\vol 2
\yr 1961
\pages 728--733
\endref

\ref
\noOf{41}
\by Vyssotsky, V.\ A., Gordon, S.\ B., Frisch, H.\ L., \jmh 
\paper
Critical percolation probabilities (bond problem)
\jour\TPR
\vol 123
\yr 1961
\pages 1566--1567
\endref

\ref
\noOf{n18}
\by Frisch, H.\ L., \jmh, Sonnenblick, E., Vyssotsky, V.\ A.
\paper Critical probabilities: site problem
\jour \TPR
\vol 124
\yr 1961
\pages 1021--1022
\endref

\ref
\noOf{n19}
\by\jmh
\paper On Steiner's network problem
\jour Mathematika
\vol 8
\yr 1961
\pages 131--132
\endref

\ref
\noOf{42}
\by \jmh
\paper
On the statistical loss of long-period comets from the solar system. II.
\jour Proceedings of the Fourth Berkeley Symposium on Mathematics,  
Statistics and Probability
\vol III
\yr 1961
\pages 17--78
\endref

\ref
\noOf{44}
\by \jmh
\paper
On the dynamical disequilibrium of individual particles
\jour  Proceedings of the Fourth Berkeley Symposium on Mathematics, Statistics
and Probability
\vol III
\yr 1961
\pages 79--85
\endref

\ref
\noOf{43}
\by \jmh
\paper
On the rate of convergence to the connective constant of the hypercubical
lattice
\jour The Quarterly Journal of Mathematics. Oxford 
\vol 12 
\yr 1961
\pages 250--256
\endref

\nref
\noOf{n20}
\by \jmh, Levine, H.
\paper Planning for the distant future
\jour The Times Educational Supplement
\yr 15 September 1961
\page 293
\endref

\ref
\noOf{n21}
\by Cranshaw, T.\ E., \jmh
\paper Counting statistics
\jour Encyclopaedic Dictionary of Physics
\vol 2
\yr 1962
\pages 89--108
\endref

\ref
\noOf{n22}
\by Frisch, H.\ L., Gordon, S.\ B., \jmh, Vyssotsky, V.\ A.
\paper Monte Carlo solution of bond percolation processes
in various crystal lattices
\jour Bell System Technical Journal
\vol 41
\yr 1962
\pages 909--920
\endref

\ref
\noOf{n23}
\by Frisch, H.\ L. \jmh, Welsh, D.\ J.\ A.
\paper Monte Carlo estimates of percolation
probabilities for various lattices
\jour \TPR
\vol 126
\yr 1962
\pages 949--951
\endref

\ref
\noOf{45}
\by \jmh
\paper
Generalization of the fundamental theorem on sub-additive functions
\jour \PCPS 
\vol 58
\yr 1962
\pages 235--238
\endref

\ref
\noOf{46}
\by \jmh , Welsh, D.\ J.\ A.
\paper
Further results on the rate of convergence to the connective constant of the 
hypercubical lattice
\jour The Quarterly Journal of Mathematics. Oxford 
\vol 13
\yr 1962
\pages 108--110
\endref

\ref
\noOf{47}
\by \jmh
\paper
The mathematical analysis of traffic congestion
\jour Bulletin de l'Institut International de Statistique 
\vol 39 
\yr 1962
\pages 89--108
\endref

\ref
\noOf{n24}
\by\jmh
\paper Monte Carlo methods
\inbook Proceedings of the 7th Conference on the Design of Experiments
in Army Research
\publ Development and Testing, U.\ S.\ Army Research Office
\yr 1962
\pages 17--26
\endref

\ref
\noOf{n25}
\by\jmh
\paper A Monte Carlo solution of percolation in the cubic lattice
\jour Meth.\ Comput.\ Phys.
\vol 1
\yr 1963
\pages 281--298
\endref

\ref
\noOf{48}
\by Frisch, H.\ L., \jmh 
\paper
Percolation processes and related topics
\jour Journal of the Society for Industrial and Applied Mathematics
\vol 11
\yr 1963
\pages 894--918
\endref

\ref
\noOf{49}
\by \jmh , Walters, R.\ S.
\paper
Percolation and fractional branching processes
\jour \JSIAM
\vol 11
\yr 1963
\pages 831--839
\endref

\ref
\noOf{50}
\by \jmh
\paper
Long-chain polymers and self-avoiding random walks
\jour  Sankhya, Series A
\vol 25
\yr 1963
\pages 29--38
\endref

\ref
\noOf{51}
\by \jmh
\paper
Long-chain polymers and self-avoiding random walks. II.
\jour Sankhya, Series A
\vol 25 
\yr 1963
\pages 269--272
\endref

\nref
\noOf{53}
\by \jmh, Handscomb, D.\ C.
\book
Monte Carlo Methods
\publ Methuen 
\publaddr London
\yr 1964
\endref

\ref
\noOf{52}
\by \jmh , Welsh, D.\ J.\ A.
\paper
First-passage percolation, subadditive processes, stochastic networks, and
generalized renewal theory
\inbook Bernoulli, Bayes, Laplace Anniversary Volume
\eds Neyman, J., LeCam, L.
\yr 1965
\pages 61--110
\endref

\ref
\noOf{54}
\by \jmh
\paper
Subadditive functional expectations
\jour Theory of Probability and its Applications (Russian)
\vol 11
\yr 1966
\pages 352--354
\moreref 311--313 (English)
\endref

\ref
\noOf{55}
\by \jmh
\paper
Existence theorems and Monte Carlo methods for the monomer--dimer problem
\inbook  Research Papers in Statistics (Festschrift for J.\ Neyman)
\publ John Wiley
\publaddr London
\yr 1966
\pages 125--146 
\endref

\ref
\noOf{56}
\by \jmh
\paper
First-passage percolation
\jour \JRSSB
\vol 28
\yr 1966
\pages 491--496
\endref

\ref
\noOf{57}
\by \jmh , Mallows, C.\ L., Handscomb, D.\ C.
\paper
Recent publications and presentations: Monte Carlo methods
\jour The American Mathematical Monthly 
\vol 73 
\yr 1966 
\pages 685
\endref

\ref
\noOf{58}
\by Bingham, N.\ H., \jmh
\paper
On a conjecture of Rademacher, Dickson, and Plotkin
\jour Journal of  Combinatorial Theory
\vol 3 
\yr 1967
\pages 182--190
\endref

\ref
\noOf{59}
\by \jmh
\paper
Harnesses
\jour  Proceedings of the Fifth Berkeley Symposium on Mathematical 
Statistics and Probability 
\vol III
\yr 1967
\pages pp. 89--117
\endref

\nref
\noOf{n28}
\by\jmh
\paper Industry and education: prospects and responsibilities
in mathematics in South Africa
\jour Kwart.\ Tyd.\ Wisk.\ Wetenskap.
\vol 5
\yr 1967
\pages 11--17
\endref

\nref
\noOf{60}
\by \jmh , Handscomb, D.\ C.
\paper
Les m\'ethodes de Monte-Carlo
\inbook Traduit de l'anglais par F. Rostand. Monographies Dunod, No. 65
\publ Dunod
\publaddr Paris 
\yr 1967
\pages 229 
\endref

\nref
\noOf{new1}
\by\jmh
\paper On the enfeeblement of mathematical
skills by `Modern Mathematics' and by similar soft
intellectual trash in schools and universities
\jour \BIMA
\vol 4
\pages 66--85
\yr 1968
\endref

\ref
\noOf{61}
\by \jmh
\paper
An improved lower bound for the multidimensional dimer problem
\jour \PCPS
\vol 64 
\yr 1968
\pages 455--463
\endref

\ref
\noOf{n29}
\by Feuerverger, A., \jmh, Izenman, A., Makani, K.
\paper Negative finding for the three-dimensional
dimer problem
\jour \JMP
\vol 10
\yr 1969
\pages 443--446
\endref

\ref
\noOf{n45}
\by \jmh
\paper
Sequences of absolute differences
\jour SIAM Review
\vol 11
\yr 1969
\pages 73--74
\endref

\ref
\noOf{n30}
\by \jmh
\paper Calculation of lattice statistics
\inbook Proceedings of the 2nd Conference on Computational
Physics
\publ
Institute of Physics and Physical Society
\publaddr London
\yr 1970
\pages 1--8
\endref

\ref
\noOf{62}
\by \jmh , Menon, V.\ V.
\paper
A lower bound for the monomer--dimer problem
\jour Journal of the Institute of Mathematics and its Applications 
\vol 6
\yr 1970
\pages 341--364
\endref

\ref
\noOf{HamC}
\by\jmh, Clifford, P.
\paper Markov fields on finite graphs and lattices
\yr 1971
\paperinfo unpublished
\endref

\nref
\noOf{n31}
\by\jmh
\paper No matter, never mind!
paperinfo Lecture deliverd at the dedication of Evans Hall,
U.\ C.\ Berkeley
\jour \BIMA
\vol 7
\yr 1971
\pages 358--364
\endref

\ref
\noOf{63}
\by \jmh
\paper
A few seedlings of research
\jour Proceedings of the Sixth Berkeley Symposium on Mathematical Statistics 
and Probability
\vol I
\yr 1972
\pages 345--394
\endref

\nref
\noOf{n32}
\by\jmh
\paper Symposium on teaching of mathematics in
schools in relation to the teaching of physics (Eton College,
2 October 1971): impression of
the meeting
\jour \BIMA
\vol 8
\yr 1972
\pages 39--40
\endref

\ref
\noOf{n34}
\by\jmh
\paper Stochastic models for the distribution of particles in space
\jour \AAP\ (Supplement)
\yr 1972
\pages 47--68
\endref

\ref
\noOf{new50}
\by\jmh
\paper Maxims for manipulators
\jour\BIMA
\vol 9
\yr 1973
\pages 276--280
\endref

\nref
\noOf{n35}
\by\jmh
\paper How is research done?
\jour\BIMA
\vol 9
\yr 1973
\pages 214--215
\endref

\ref
\noOf{n33}
\by\jmh
\paper Contribution to discussion on subadditive
ergodic theory
\jour \AP
\vol 1 
\yr 1973
\pages 905--909
\endref

\nref
\noOf{n36}
\by\jmh
\paper Statistical tools
\jour The Statistician
\vol 23
\yr 1974
\pages 89--106
\endref

\ref
\noOf{n38}
\by Bell, G.\ M., Churchhouse, R.\ F., Goodwin, E.\ T., \jmh, Taylor, R.\ S.
\paper Proof of a conjecture of Worster
\jour\BIMA
\vol 10
\yr 1974
\pages 128--129
\endref

\nref
\noOf{new3}
\by\jmh
\paper Poking about for the vital juices of mathematical research
\jour\BIMA
\vol 10
\yr 1974
\pages 235--247
\endref

\nref
\noOf{new4}
\by\jmh
\paper Some thoughts occasioned by an undergraduate mathematics society
\jour\BIMA
\yr 1974
\vol 10
\pages 306--311
\endref

\ref
\noOf{new5}
\by\jmh
\paper An isoperimetric problem
\jour\BIMA
\vol 10
\pages 439--441
\yr 1974
\endref

\ref
\noOf{new6}
\by\jmh
\paper A rather difficult O-level problem
\jour\BIMA
\vol 10
\pages 441--443
\yr 1974
\endref

\nref
\noOf{n39}
\by\jmh
\paper The technology of thought
\inbook
The Heritage of Copernicus
\yr 1974
\ed J.\ Neyman
\publ MIT Press
\endref

\ref
\noOf{65}
\by \jmh
\paper
Postulates for subadditive processes
\jour Annals of Probability 
\vol 2 
\yr 1974 
\pages 652--680
\endref

\ref
\noOf{66a}
\by \jmh
\paper
Some speculations on a sense of nicely calculated chances
\jour SIAM Review
\vol 16 
\yr 1974
\pages 237--255
\endref

\ref
\noOf{66}
\by \jmh, Grimmett, G.\ R.
\paper
Maximal solutions of the generalized subadditive inequality
\inbook Stochastic Geometry (A tribute to the memory of Rollo Davidson)
\publ Wiley
\publaddr London
\yr 1974
\pages 270--284
\endref

\nref
\noOf{67}
\by Besag, Julian
\paper
Spatial interaction and the statistical analysis of lattice systems. With 
discussion by D. R. Cox, A. G. Hawkes, P. Clifford, P. Whittle, K. Ord, 
R. Mead, J. M. Hammersley, and M. S. Bartlett and with a reply by the author. 
\jour \JRSSB 
\vol  36 
\yr 1974
\pages 192--236
\endref

\nref
\noOf{new2}
\by\jmh
\paper Lehrs\"atze and Leers\"atze di Polenta e Segu
\jour\BIMA
\vol 11
\pages 117--121
\yr 1975
\endref

\ref
\noOf{n41}
\by\jmh
\paper Some general reflections on statistical practice
\inbook Festschrift for Professor Linder
\endref

\ref
\noOf{n40}
\by\jmh
\paper The wide open spaces
\jour The Statistician
\vol 24
\yr 1975
\pages 159--160
\endref

\ref
\noOf{68}
\by \jmh
\paper
Rumination on infinite Markov systems
\inbook  Perspectives in Probability and 
Statistics (Papers in honour of M. S. Bartlett)
\publ Applied Probability Trust
\publaddr Sheffield
\yr 1975
\pages 195--200
\endref

\ref
\noOf{69}
\by \jmh , Lewis, J.\ W.\ E., Rowlinson, J.\ S.
\paper
Relationships between the multinomial and Poisson models of stochastic 
processes, and between the canonical and grand canonical ensembles in 
statistical mechanics, with illustrations and Monte Carlo methods for the 
penetrable sphere model of liquid--vapour equilibrium
\jour Sankhya, Series A
\vol 37
\yr 1975 
\pages 457--491
\endref

\ref
\noOf{70}
\by \jmh
\paper
The design of future computing machinery for functional integration
\inbook Functional Integration and its Applications
\publ Clarendon Press
\publaddr Oxford
\yr 1975
\pages 83--86
\endref

\ref
\noOf{71}
\by \jmh
\paper
A generalization of McDiarmid's theorem for mixed Bernoulli percolation
\jour Mathematical Proceedings of the Cambridge Philosophical Society
\vol 88
\yr 1980 
\pages 167--170
\endref

\ref
\noOf{72}
\by \jmh
\paper
Biological growth and spread
\inbook  Lecture Notes in Biomathematics
\vol 38
\yr 1980
\publ Springer
\publaddr Berlin
\pages 484--494
\endref

\ref
\noOf{new80}
\by \jmh, Welsh, D.\ J.\ A.
\paper Percolation theory and its ramifications
\jour Contemporary Physics
\vol 21
\yr 1980
\pages 593--505
\endref

\ref
\noOf{73}
\by \jmh
\paper
Critical phenomena in semi-infinite systems. Essays in statistical science
\jour \JAP 
\vol 19A
\yr 1982
\pages 327--331
\endref

\ref
\noOf{74}
\by \jmh , Torrie, G.\ M., Whittington, S.\ G.
\paper
Self-avoiding walks interacting with a surface
\jour Journal of Physics. A. Mathematical and General 
\vol 15 
\yr 1982 
\pages 539--571
\endref

\nref
\noOf{75}
\by \jmh
\paper
Obituary: J. Neyman, 1894--1981
\jour \JRSSA 
\vol 145 
\yr 1982 
\pages 523--524
\endref

\ref
\noOf{76}
\by \jmh
\paper
Oxford commemoration ball
\inbook Probability, Statistics and Analysis
\bookinfo London Mathematical Society Lecture Note Series 
\vol 79
\publ Cambridge University Press
\publaddr Cambridge
\yr 1983
\pages 112--142
\endref

\ref
\noOf{77}
\by \jmh
\paper
The friendship theorem and the love problem
\inbook Surveys in Combinatorics
\bookinfo London Mathematical Society Lecture Note Series 
\vol 82
\publ Cambridge University Press
\publaddr Cambridge
\yr 1983
\pages 31--54
\endref

\ref
\noOf{78}
\by \jmh
\paper
Origins of percolation theory
\inbook Percolation structures and processes
\bookinfo Annals of the Israel Physical Society 
\vol 5
\publ Hilger
\publaddr Bristol
\yr 1983
\pages 47--57
\endref

\ref
\noOf{79}
\by \jmh , Mazzarino, G.
\paper
Markov fields, correlated percolation, and the Ising model
\inbook The Mathematics and Physics of Disordered Media  
\bookinfo Lecture Notes in Mathematics
\vol 1035
\publ Springer
\publaddr Berlin
\yr 1983
\pages 210--245
\endref

\ref
\noOf{80}
\by \jmh
\paper
Functional roots and indicial semigroups
\jour Bulletin of the Institute of Mathematics and its Applications
\vol 19
\yr 1983 
\pages 194--196
\endref

\ref
\noOf{81}
\by \jmh , Whittington, S.\ G.
\paper
Self-avoiding walks in wedges
\jour Journal of Physics. A. Mathematical and General 
\vol 18 
\yr 1985
\pages 101--111
\endref

\ref
\noOf{82}
\by \jmh
\paper
Mesoadditive processes and the specific conductivity of lattices
\inbook A Celebration of Applied Probability
\bookinfo \JAP 
\vol 25A
\yr 1988
\pages 347--358
\endref

\nref
\noOf{82a}
\by \jmh 
\paper
Room to wriggle
\jour Bulletin of the Institute of Mathematics and its Applications 
\vol 24 
\yr 1988 
\pages 65--72
\endref

\ref
\noOf{83}
\by \jmh , Mazzarino, G.
\paper
A differential equation connected with the dendritic growth of crystals
\jour IMA Journal of Applied Mathematics 
\vol 42 
\yr 1989 
\pages 43--75
\endref

\ref
\noOf{84}
\by \jmh , Mazzarino, G.
\paper
Computational aspects of some autonomous differential equations
\jour Proceedings of the Royal Society, Series A 
\vol 424 
\yr 1989
\pages 19--37
\endref

\nref
\noOf{85}
\book 
Disorder In Physical Systems. A volume in honour of John M. Hammersley on the
occasion of his 70th birthday
\eds G.\ R. Grimmett and D.\ J.\ A. Welsh
\publ Clarendon Press
\publaddr New York
\yr 1990
\endref

\nref
\noOf{86}
\by Kendall, David
\paper 
Speech proposing the toast to John Hammersley---1 October 1987
\inbook Disorder in physical systems
\eds G.\ R. Grimmett and D.\ J.\ A. Welsh
\publ  Oxford University Press
\publaddr New York
\yr 1990
\pages 1--3
\endref

\ref
\noOf{87}
\by \jmh
\paper Self-avoiding walks
\inbook Current problems in statistical mechanics 
\bookinfo Physica A 
\vol 177
\yr 1991
\pages 51--57
\moreref Corrigendum  vol.\ 183 (1992), 574--578 
\endref

\nref
\noOf{88}
\by \jmh
\paper
Corrigendum: ``Self-avoiding walks'' [Phys. A 177 (1991), no. 1-3, 51--57]
\jour Physica A. Statistical and Theoretical Physics 
\vol 183, no.4
\yr 1992 
\pages 574--578
\endref

\ref
\noOf{89}
\by \jmh , Mazzarino, G.
\paper
Properties of large Eden clusters in the plane
\jour Combinatorics, Probability and Computing  
\vol 3
\yr 1994
\pages 471--505
\endref

\ref
\noOf{90}
\by \jmh
\paper
Fractal dynamics of Eden clusters
\inbook Probability, Statistics and Optimisation
\publ Wiley
\publaddr Chichester
\yr 1994
\pages 79--87
\endref
 
\subrefhead{Publications on other topics}

\ref
\noOf{n13}
\book Proceedings of the Oxford Mathematical Conference
for Schoolteachers and Industrialists
\ed \jmh
\yr 1957
\publ Times Publishing Company
\publaddr London
\endref

\ref
\noOf{n14}
\by\jmh
\paper The value of mathematics and its teachers
\jour ibid
\yr 1957
\endref

\ref
\noOf{n17}
\by Coulson, C.\ A., \jmh
\paper The bottleneck in British science and technology
\jour New Scientist
\vol 10
\yr 1961
\pages 499--500
\endref

\ref
\noOf{n20}
\by \jmh, Levine, H.
\paper Planning for the distant future
\jour The Times Educatonal Supplement
\yr 15 September 1961
\page 293
\endref

\ref
\noOf{n28}
\by\jmh
\paper Industry and education: prospects and responsibilities
in mathematics in South Africa
\jour Kwart.\ Tyd.\ Wisk,.\ Wetenskap.
\vol 5
\yr 1967
\pages 11--17
\endref

\nref
\noOf{60}
\by \jmh , Handscomb, D.\ C.
\paper
Les m\'ethodes de Monte-Carlo
\inbook Traduit de l'anglais par F. Rostand. Monographies Dunod, No. 65
\publ Dunod
\publaddr Paris
\yr 1967
\pages 229
\endref

\ref
\noOf{new1}
\by\jmh
\paper On the enfeeblement of mathematical
skills by `Modern Mathematics' and by similar soft
intellectual trash in schools and universities
\jour \BIMA
\vol 4
\pages 66--85
\yr 1968
\endref

\ref
\noOf{n31}
\by\jmh
\paper No matter, never mind!
\jour \BIMA
\vol 7
\yr 1971
\pages 358--364
\endref

\ref
\noOf{n32}
\by\jmh
\paper Symposium on teaching of mathematics in
schools in relation to the teaching of physics (Eton College,
2 October 1971): impression of
the meeting
\jour \BIMA
\vol 8
\yr 1972
\pages 39--40
\endref

\ref
\noOf{n35}
\by\jmh
\paper How is research done?
\jour\BIMA
\vol 9
\yr 1973
\pages 214--215
\endref

\ref
\noOf{n65}
\by\jmh
\paper Modern mathematics, the great debate:
Motion proposing that this house deplores the enthusiastic
teaching of modern mathematics, particularly in schools
\jour\BIMA
\vol 9
\yr 1973
\pages 238--241
\endref

\ref
\noOf{new3}
\by\jmh
\paper Poking about for the vital juices of mathematical research
\jour\BIMA
\vol 10
\yr 1974
\pages 235--247
\endref

\ref
\noOf{new4}
\by\jmh
\paper Some thoughts occasioned by an undergraduate mathematics society
\jour\BIMA
\yr 1974
\vol 10
\pages 306--311
\endref

\ref
\noOf{n36}
\by\jmh
\paper Statistical tools
\jour The Statistician
\vol 23
\yr 1974
\pages 89--106
\endref

\ref
\noOf{n39}
\by\jmh
\paper The technology of thought
\inbook
The Heritage of Copernicus
\yr 1974
\ed J.\ Neyman
\publ MIT Press
\endref

\ref
\noOf{new2}
\by\jmh
\paper Lehrs\"atze and Leers\"atze di Polenta e Segu
\jour\BIMA
\vol 11
\pages 117--121
\yr 1975
\endref

\ref
\noOf{new51}
\by\jmh
\paper Sweet nothing
\jour \BIMA
\vol 14
\yr 1978
\pages 146--147
\endref

\ref
\noOf{75}
\by \jmh
\paper
Obituary: J. Neyman, 1894--1981
\jour \JRSSA
\vol 145
\yr 1982
\pages 523--524
\endref

\ref
\noOf{new52}
\by\jmh
\paper The teaching of combinatorial analysis
\jour\BIMA
\vol 19
\yr 1983
\pages 50--52
\endref

\ref
\noOf{new53}
\by\jmh
\paper Probability and arithmetic in science
\jour\BIMA
\vol 21
\yr 1985
\pages 114--120
\endref

\ref
\noOf{82a}
\by \jmh
\paper
Room to wriggle
\jour Bulletin of the Institute of Mathematics and its Applications
\vol 24
\yr 1988
\pages 65--72
\endref

\endRefs

\bye